\documentclass[doublespacing]{elsart}

\usepackage{graphicx}
\usepackage{amssymb}

\usepackage{amsmath}
\usepackage{amsfonts}

\newcommand{\Naturales}[0] {\mathbb{N}}

\newcommand{\dis}[0] {\displaystyle}

\begin{document}

\begin{frontmatter}

\title{A curious arithmetic of fractal dimension for polyadic Cantor sets}

\author{Francisco R. Villatoro\corauthref{COR}}

 \address{Departamento de Lenguajes y Ciencias de la
 Computaci\'on,\\
  Universidad de M\'alaga, E-29071 M\'alaga, Spain}

 \corauth[COR]{Corresponding author: villa@lcc.uma.es}

\begin{abstract}

Fractal sets, by definition, are non-differentiable, however their
dimension can be continuous, differentiable, and arithmetically
manipulable as function of their construction parameters. A new
arithmetic for fractal dimension of polyadic Cantor sets is
introduced by means of properly defining operators for the addition,
subtraction, multiplication, and division. The new operators have
the usual properties of the corresponding operations with real
numbers. The combination of an infinitesimal change of fractal
dimension with these arithmetic operators allows the manipulation of
fractal dimension with the tools of calculus.
\end{abstract}


\end{frontmatter}

%%%%%%%%%%%%%%%%%%%%%%%%%%%%%%%%%%%%%%%%
%%%%%%%%%%%%%%%%%%%%%%%%%%%%%%%%%%%%%%%%

\section{Introduction}
\label{Intro}

This paper is motivated by the recent introduction of a series of
arithmetic operators for fractal dimension of
hyperhelices~\cite{Fletcher2001,Fletcher2004} by Carlos D.
Toledo-Su{\'a}rez~\cite{Carlos2007}. Let us briefly recall that a
hyperhelix of order $N$ is defined to be a self-similar object
consisting of a thin elastic rod wound into a helix, which is itself
wound into a larger helix, until this process has been repeated $N$
times. A hyperhelix fractal results when $N$ tends to infinity. The
Hausdorff-Besicovitch dimension of a hyperhelix fractal can be any
desired number above unity. However, the interpretation of
hyperhelices with fractal dimension larger than three, which
unavoidably intersect themselves, is not obvious. In fact, this
author is attempted to point out that geometrical considerations
forbids hyperhelix fractal dimensions larger than three. This is the
main inconvenient with the fractal dimension arithmetics introduced
by Toledo-Su{\'a}rez~\cite{Carlos2007}. Could such drawbacks be
circumvented by using another kind of fractal?

This short paper introduces an algebraic structure for the fractal
dimension of polyadic Cantor fractals, resulting from iterative
replacement of segments in the unit interval by scaled versions of
them, whose Hausdorff-Besicovitch dimension can be any real number
in the interval $[0,1]$. The new arithmetic preserves this property
avoiding the appearance of non-geometrically feasible fractal
dimensions larger than unity. The contents of this paper are as
follows. The next section recalls the characteristics of
symmetrical, polyadic, Cantor fractal sets. Section~3, presents the
new arithmetical operators for the addition, subtraction, product,
and division of the fractal dimension, summarizing their main
properties. Finally, the last section is devoted to the main
conclusions.

\section{Polyadic Cantor fractals}
\label{Superlattices}

\begin{figure}%%[ht]
 \begin{center}
  \includegraphics[width=12cm]{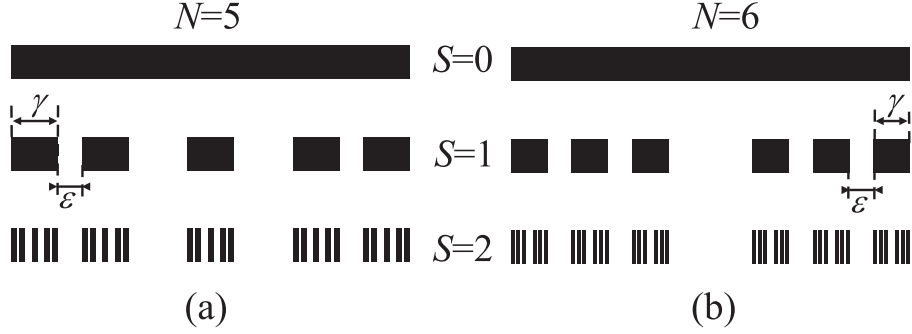}
 \end{center}
  \caption{First two stages ($S=1$ and $S=2$) of
  polyadic Cantor fractal sets with $N=5$ (a) and $N=6$ (b),
  showing the definition
  of both the scale factor ($\gamma$) and the lacunarity parameter
  ($\varepsilon$).}
 \label{PolyadicCantorSet}
\end{figure}

Polyadic, or generalized, Cantor sets are defined as
follows~\cite{Mandelbrot,Villatoro2006}. The first step ($S=0$) is
to take a segment of unit length, referred to as the initiator. In
the next step, $S=1$, the segment is replaced by $N$ non-overlapping
copies of the initiator, each one scaled by a factor $\gamma<1$. For
$N$ odd, as shown in Fig.~\ref{PolyadicCantorSet}(a), one copy lies
exactly centered in the interval, and the remaining ones are
distributed such that $\lfloor N/2 \rfloor$ copies are placed
completely to the left of the interval and the remaining $\lfloor
N/2 \rfloor$ copies are placed wholly to its right, where $\lfloor
N/2 \rfloor$ is the greatest integer less than or equal to $N/2$;
additionally, each copy among these $N-1$ ones is separated by a
fixed distance, let us say $\varepsilon$.  For $N$ even, as shown in
Fig.~\ref{PolyadicCantorSet}(b), one half of the copies is placed
completely to the left of the interval and the other half completely
to its right, with each copy separated by $\varepsilon$. At the
following construction stages of the the fractal set, $S=2, 3,
\ldots$, the generation process is repeated over and over again for
each segment in the previous stage. Strictly speaking, the Cantor
set is the limit of this procedure as $S\rightarrow\infty$, which is
composed of geometric points distributed so that each point lies
arbitrarily close of other points of the set, being the $S$-th stage
Cantor set usually referred to as a pre-fractal or physical fractal.

Symmetrical polyadic Cantor fractals are characterized by three
parameters, the number of self-similar copies $N$, the scaling
factor $\gamma$, and the width of the outermost gap at the first
stage, $\varepsilon$, here on referred to as the lacunarity
parameter, as in most of the previous technical papers dealing with
polyadic Cantor fractals in engineering
applications~\cite{JaggardJOSA,MonsoriuOE,Villatoro2008}. The
similarity dimension of all polyadic Cantor fractals is
$D=\ln(N)/\ln(\gamma^{-1})=-\ln(N)/\ln(\gamma)$, which is
independent of the lacunarity parameter.

The three parameters of a polyadic Cantor set must satisfy certain
constraints in order to avoid overlapping between the copies of the
initiator. Let us recall such almost obvious results for the sake of
completeness. On the one hand, the maximum value of the scaling
factor depends on the value of $N$, such that $0<\gamma
<\gamma_{\mbox{\scriptsize max}}= 1/N$. On the other hand, for each
$N$ and $\gamma$, there are two extreme values for $\varepsilon$.
The first one is $\varepsilon_{\mbox{\scriptsize min}}=0$, for which
the highest lacunar fractal is obtained, i.e., that with the largest
possible gap. For $N$ even, the central gap has a width equal to
$1-N\,\gamma$, and, for $N$ odd, both large gaps surrounding the
central well have a width of $(1-N\,\gamma)/2$. The other extreme
value is
\[
 \varepsilon_{\mbox{\scriptsize max}}=
 \left\{ \begin{array}{ll}
 \dis \frac{1-N\,\gamma}{N-2},\qquad & \mbox{even } N, \\
  \\
 \dis \frac{1-N\,\gamma}{N-3},\qquad & \mbox{odd } N,
 \end{array}\right.
\]
where for even (odd) $N$ two (three) wells are joined together in
the center, without any gap in the central region. The width of the
$N-2$ gaps in this case is equal to $\varepsilon_{\mbox{\scriptsize
max}}$. Thus, the corresponding lacunarity is lower than that for
$\varepsilon=0$, but not the smallest one, which is obtained for the
most regular distribution, where the gaps and wells have the same
width at the first stage ($S=1$) given by
\[
 \varepsilon_{\mbox{\scriptsize reg}}= \frac{1-N\,\gamma}{N-1}.
\]
Note that $0 < \varepsilon_{\mbox{\scriptsize reg}} <
\varepsilon_{\mbox{\scriptsize max}}$.

\section{Arithmetic of fractal dimension}

Let us take two $N$-adic Cantor fractals $A$ and $B$ with scaling
factors (fractal dimensions) $\gamma_A$ ($D_A$) and $\gamma_B$
($D_B$) , respectively. Their fractal dimensions can be combined by
means of the following four operators in order to obtain a new
$N$-adic Cantor fractal $C$ with scaling factor $\gamma_C$ and
fractal dimension $D_C$. Note that the independency of the
similarity dimension on the lacunarity parameter avoids the need of
its incorporation into the dimension arithmetic which follows.

\subsection{Addition ($\oplus$)}

Let us define the $N$-adic Cantor fractal $C=A\oplus B$ by means of
$\gamma_C=\gamma_A\,\gamma_B$, resulting in
\[
 \frac{1}{D_C} = \frac{1}{D_A} + \frac{1}{D_B},
 \qquad
 D_C = D_A\oplus D_B = \frac{D_A\,D_B}{D_A+D_B},
\]
where, by abuse of notation, the operator $\oplus$ has also been
applied to the fractal dimensions of the fractals $A$ and $B$. The
present definition of the operator $\oplus$ is consistent, since for
$0<\gamma_A,\gamma_B<1/N$, then $0<\gamma_C<1/N$, i.e., for
$0<D_A,D_B<1$, then $0<D_C<1$.

The operator $\oplus$ is commutative ($D_A\oplus D_B = D_B\oplus
D_A$) and associative
\[
 D_A\oplus (D_B \oplus D_C)
 =
 \frac{1}{\frac{1}{D_A}+\frac{1}{D_B}+\frac{1}{D_C}}
 =
 (D_A\oplus D_B) \oplus D_C.
\]

The void set, i.e., the $N$-adic Cantor fractal $Z$ with
$\gamma_Z=0$ and $D_Z=0$, is the absorbing element of $\oplus$,
since $D_Z=D_A\oplus D_Z$ for all $A$.

There is no identity element $Y$ for $\oplus$, such that
$D_A=D_A\oplus D_Y$, since the only possibility is to take
$\gamma_Y=1$ and $D_Y=+\infty$, which it is not a geometrically
valid $N$-adic Cantor fractal.

Let $U$ be the unit segment, with $D_U=1$ and $\gamma_U=1/N$. Hence
$D_U\oplus D_U= 1/2$ and $D_A\oplus D_U=D_A/(1+D_A)$. Note also that
$D_A\oplus D_A=D_A/2$.

\subsection{Subtraction ($\ominus$)}

Let us define the $N$-adic Cantor fractal $C=A\ominus B$, by means
of $\gamma_C=\gamma_A/ \gamma_B$, resulting in
\[
 \frac{1}{D_C} = \frac{1}{D_A} - \frac{1}{D_B},
 \qquad
 D_C = D_A\ominus D_B = \frac{D_A\,D_B}{D_B-D_A}.
\]
This definition is consistent only for $\gamma_A<\gamma_B/N$, i.e.,
for
$$D_A<D_B\oplus D_U=\frac{D_B}{1+D_B}<\frac{1}{2},$$ where
$0<D_B<1$. Under these assumptions, the operator $\ominus$ is
compatible with $\oplus$ since
\[
 (D_A \ominus D_B) \oplus D_B = D_A.
\]
However, let us note that both $(D_A \oplus D_B) \ominus D_B$, and
$D_A\ominus D_A$ are not properly defined, since $D_A \oplus D_B
> D_B \oplus D_U$, and $D_A>D_A\oplus D_U$, respectively.

The operator $\ominus$ is neither commutative ($D_A\ominus D_B \ne
D_B\ominus D_A$) nor associative
\[
 D_A\ominus (D_B \ominus D_C)
 =
 \frac{1}{\frac{1}{D_A}-\frac{1}{D_B}+\frac{1}{D_C}}
 \ne
 \frac{1}{\frac{1}{D_A}-\frac{1}{D_B}-\frac{1}{D_C}}
 =
 (D_A\ominus D_B) \ominus D_C.
\]

The void set $Z$, with $\gamma_Z=0$ and $D_Z=0$, is also the
absorbing element of $\ominus$, since $D_Z=D_A\ominus D_Z$ for all
$A$.

\subsection{Product ($\otimes$)}

Let us define the $N$-adic Cantor fractal $C=A\otimes B$ by means of
$\gamma_C=\gamma_A^{1/D_B}$, resulting in
\[
 \frac{1}{D_C} = -\frac{1}{D_B}\frac{\ln (\gamma_A)}{\ln (N)} =
 \frac{1}{D_A\,D_B},
 \qquad
 D_C = D_A\otimes D_B = D_A\,D_B.
\]
The correctness of this definition is straightforward since for
$0<\gamma_A,\gamma_B<1/N$, $0<\gamma_C<1/N$; similarly,
$0<D_A,D_B<1$ results in $0<D_C<1$.

The operator $\otimes$ is commutative ($D_A\otimes D_B = D_B\otimes
D_A$) and associative
\[
 D_A\otimes (D_B \otimes D_C)
 =
 D_A\,D_B\,D_C
 =
 (D_A\otimes D_B) \otimes D_C.
\]
It is also distributive with respect to the addition, since
\[
 D_A\otimes(D_B\oplus D_C) = \frac{D_A\,D_B\,D_C}{D_B+D_C}
 = (D_A\otimes D_B)\oplus (D_A\otimes D_C),
\]
and also with respect to the subtraction, as
\[
 D_A\otimes(D_B\ominus D_C) = \frac{D_A\,D_B\,D_C}{D_C-D_B}
 = (D_A\otimes D_B)\ominus (D_A\otimes D_C).
\]

Note that the operator $\oplus$ is not distributive with respect to
the product, since
\[
 D_A\oplus(D_B\otimes D_C) = \frac{D_A\,D_B\,D_C}{D_A+D_B\,D_C}
 \ne
 \frac{D_A^2\,D_B\,D_C}{2\,D_A+D_B\,D_C}
 = (D_A\oplus D_B)\otimes (D_A\oplus D_C).
\]
The operator $\ominus$ is also not distributive with respect to the
product.

The void set $Z$ is the absorbing element of $\otimes$, since
$D_Z=D_A\otimes D_Z$ for all $A$, and the unit segment $U$ is the
unit element of the product, since $D_A\otimes D_U= D_A$.

\subsection{Division ($\oslash$)}

Let us define the $N$-adic Cantor fractal $C=A\oslash B$ by means of
$\gamma_C=\gamma_A^{D_B}=\ln(\gamma_B)/\ln(\gamma_A)$, resulting in
\[
 \frac{1}{D_C} = -{D_B}\frac{\ln (\gamma_A)}{\ln (N)} =
 \frac{1}{D_B\,D_A},
 \qquad
 D_C = D_A\oslash D_B = \frac{D_A}{D_B}.
\]
The correctness of this definition requires that
$0<\gamma_A<\gamma_B<1/N$, i.e., $0<D_A<D_B<1$. The division
operator is compatible with the product, since
\[
 (D_A \oslash D_B)\otimes D_B =  (D_A \otimes D_B)\otimes D_B = D_A,
\]
and with the unit element $D_A\oslash D_U=D_A$.

The operator $\oslash$ is neither commutative ($D_A\oslash D_B \ne
D_B\oslash D_A$) nor associative
\[
 D_A\oslash (D_B \oslash D_C)
 =
 \frac{D_A\,D_C}{D_B}
 \ne
 \frac{D_A}{D_B\,D_C}
 =
 (D_A\oslash D_B) \oslash D_C.
\]
It is also distributive with respect to the product, since
\[
 D_A\oslash(D_B\otimes D_C) = \frac{D_A}{D_B\,D_C}
 = (D_A\oslash D_B)\otimes (D_A\oslash D_C),
\]
but is not distributive with respect to either the addition nor the
subtraction.

\begin{figure}%%[ht]
 \begin{center}
  \includegraphics[width=12cm]{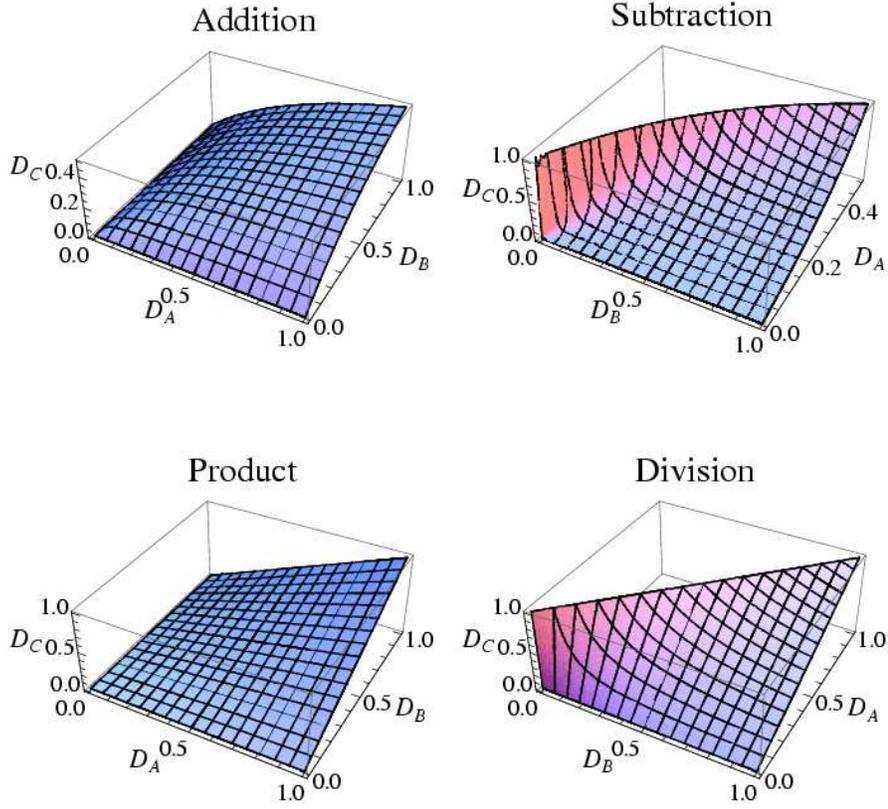}
 \end{center}
  \caption{Plots of the result of operations $D_A\oplus D_B$ (Addition),
  $D_A\ominus D_B$ (Subtraction), $D_A\otimes D_B$ (Product), and
  $D_A\oslash D_B$ (Division).}
 \label{OperatorGraphs}
\end{figure}

\subsection{Integer power}

Let us define the integer power of a $N$-adic Cantor fractal,
$C=A^{(n)}$, where $n\ge 0$ is an integer, by means of
$\gamma_C=\gamma_A^{1/{D_A^{n-1}}}$, yielding for the fractal
dimension
\[
 D_A^{(n)} = D_A \otimes D_A \otimes \overset{(n)}{\cdots} \otimes D_A
  = (D_A)^n,
 \qquad \Naturales\ni n \ge 0,
\]
where, as usual, $D_A^{(0)} = D_U$, and  $D_A^{(1)}=D_A$. The
operator of integer power of fractal dimensions has the usual
properties of the integer power of numbers inherated from those of
the product $\otimes$.

\subsection{Infinitesimal of fractal dimension}

The development of an infinitesimal calculus for fractal dimension
can be easily accomplished by using the arithmetical operators
presented in previous sections and taking into account that an
infinitesimal change of fractal dimension ruled by changes in the
factor $\gamma$ is given by
\[
 dD = - \frac{\ln(N)}{\gamma\,\ln^2(\gamma)}\,d\gamma.
\]
This infinitesimal has only one parameter, being simpler than that
introduced by Ref.~\cite{Carlos2007} which requires the
infinitesimal variation of two parameters.

\section{Conclusions}
\label{Conclusions}

Polyadic Cantor fractals may have any arbitrary dimension in the
real interval $[0,1]$ based on the continuity of their scaling
factor $\gamma$. An algebraic structure that allows to develop an
arithmetic with operators for addition, subtraction, product, and
division of fractal dimension has been introduced. The new operators
have the usual properties of the corresponding operations with real
numbers. The introduction of an infinitesimal change of fractal
dimension combined with these arithmetical operator allow the free
manipulation of fractal dimension with the tools of calculus. Our
results show that, even taking into account that fractal sets, by
definition, are non-differentiable, it is possible to differentiate
their dimension and manipulate it arithmetically.

%\bookref{authors,title,editors,chap,pp,year}
\newcommand{\bookref}[5]{#1, {#2}, #3, #4 (#5).}

%\paperref{authors,title,journal,vol,numb,pp,year}
\newcommand{\paperref}[7]{#1, {{#3}} {\bf #4}  (#6) p. #7.}
\newcommand{\paperrefno}[8]{#1, {{#3}} {\bf #4}  (#7) p. #8.}

\end{document}